\documentclass[12pt]{article}
\usepackage{amsmath,amssymb,theorem}
\usepackage{graphicx, enumerate}
\usepackage{subfigure}
\usepackage{psfrag}
\usepackage{placeins}
\usepackage{wrapfig}
\usepackage{booktabs}
\usepackage{chngcntr}
\counterwithin{figure}{section}
\numberwithin{figure}{section}
\counterwithin{table}{section}
\usepackage{hyperref}
\hypersetup{
  colorlinks,
  citecolor=black,
  filecolor=black,
  linkcolor=black,
  urlcolor=black
}

\newcommand{\mc}{\mathcal}

\newtheorem{thm}{Theorem}[section]
\newtheorem{conj}[thm]{Conjecture}

\newtheorem{claim}[thm]{Claim}
\theorembodyfont{\rmfamily}
\def\pf{\bigskip\noindent {\emph{Proof}.}~~}

\def\dfn#1{{\sl #1}}
\def\es{\emptyset}
\def\less{\setminus}

\def\qed{ \hfill $\blacksquare$}

\begin{document}
\title{Gallai-Ramsey numbers of  $C_9$ with multiple colors}
\author{
Christian Bosse and Zi-Xia Song\thanks{Corresponding Author. Email: Zixia.Song@ucf.edu}\\
Department  of Mathematics\\
University of Central Florida\\
Orlando, FL 32816, USA\\
}
\maketitle
\begin{abstract}

We  study Ramsey-type problems in Gallai-colorings. 
Given a graph $G$ and an integer $k\ge1$, the  Gallai-Ramsey number
$gr_k(K_3,G)$ is the least positive  integer $n$ such that every $k$-coloring of  the edges of the complete graph on $n$ vertices contains either a
rainbow triangle   or a monochromatic copy of  $G$.  It turns out that $gr_k(K_3, G)$ behaves more nicely than the classical Ramsey number $r_k(G)$. However,  finding exact values of $gr_k (K_3, G)$ is  far from trivial.  In this paper, we
prove that $gr_k(K_3, C_9)= 4\cdot 2^k+1$ for all $k\ge1$.  This  new result provides partial evidence for the first open case of the Triple Odd Cycle Conjecture of  Bondy and Erd\H{o}s from 1973. Our technique relies heavily on the structural result of Gallai on edge-colorings of complete graphs without rainbow triangles. We believe the method we developed can be used to  determine the exact values of $gr_k(K_3, C_n)$ for odd integers  $n\ge11$.

\end{abstract}

\section{Introduction}
\baselineskip 18pt

In this paper, we only consider finite simple graphs. The complete graph and the cycle on $n$ vertices are denoted $K_n$ and $C_n$, respectively.  We use  $|G|$ to denote the  number of vertices of  a graph $G$. \medskip

For an integer $k\ge1$, let  $c: E(G) \to [k]$  be a $k$-edge-coloring of a complete graph $G$, where $[k]:=\{1,2, \dots, k\}$. Then $c$ is a \dfn{Gallai-coloring} of $G$ if $G$ contains no  rainbow triangle (that is, a triangle with all its edges  colored differently) under  $c$. Gallai-colorings naturally arise in several areas including in information theory~\cite{KG}, in the study of partially ordered sets, as in Gallai's original paper~\cite{gallai}, and in the study of perfect graphs~\cite{CEL}. There are now a variety of papers  which consider Ramsey-type problems in Gallai-colorings (see, e.g., \cite{chgr, c5c6,GS, exponential, upbounds}).   These works mainly focus on finding various monochromatic subgraphs in such colorings. More information on this topic  can be found in~\cite{ FGP, FMO}. \medskip

For a graph $G$ and a set $A\subseteq V(G)$, we use $G[A]$ to denote the  subgraph of $G$ obtained from $G$ by deleting all vertices in $V(G)\less A$.  A graph $H$ is an \dfn{induced subgraph} of $G$ if $H=G[A]$ for some $A\subseteq V(G)$.  Recall that the classical Ramsey number $r_k(H)$ of a graph $H$ is  the least positive  integer $n$ such that every $k$-edge-coloring of   $K_n$  contains  a monochromatic copy of  $H$.  Ramsey numbers are notoriously difficult to compute in general. In this paper, we  consider Gallai-Ramsey problems. Given a graph $H$ and an integer $k\ge1$, the  \dfn{Gallai-Ramsey number}
$gr_k(K_3, H)$ is the least positive  integer $n$ such that every $k$-edge-coloring  of $K_n$ contains either a
rainbow triangle  or a monochromatic copy of  $H$.  Clearly, $gr_k(K_3,H) \leq r_k(H)$.   
The following is a result  on the general behavior of $gr_k(K_3, H)$. 
\medskip

\begin{thm} [\cite{exponential}]
Let $H$ be a fixed graph with no isolated vertices and let $k\ge1$ be an integer. If $H$ is not bipartite, then
$gr_k (K_3, H) $ is exponential in $k$. If $H$ is bipartite, then $gr_k (K_3, H) $ is linear
in $k$.			
\end{thm}

It turns out that for some graphs $H$ (e.g., when $H=C_3$),  $gr_k(K_3, H)$ behaves nicely, while the order of magnitude  of $r_k(H)$ seems hopelessly difficult to determine.  It is worth noting that  finding exact values of $gr_k (K_3, H)$ is  far from trivial, even when $|H|$ is small.
We will utilize the following important structural result of Gallai~\cite{gallai} on Gallai-colorings of complete graphs.

\begin{thm}[\cite{gallai}]\label{Gallai}
	For any Gallai-coloring $c$ of a complete graph $G$, $V(G)$ can be partitioned into nonempty sets  $V_1, V_2, \dots, V_p$ with $p>1$ so that    at most two colors are used on the edges in $E(G)\less (E(V_1)\cup \cdots\cup  E(V_p)$ and only one color is used on the edges between any fixed pair $(V_i, V_j)$ under $c$, where $E(V_i)$ denotes the set of edges in $G[V_i]$ for all $i\in [p]$. 
\end{thm}

The partition given in Theorem~\ref{Gallai} is  a \dfn{Gallai-partition} of $G$ under  $c$.  Given a Gallai-partition $V_1, V_2, \dots, V_p$ of the complete graph $G$ under $c$, let $v_i\in V_i$ for all $i\in[p]$ and let $\mathcal{R}:=G[\{v_1, v_2, \dots, v_p\}]$. Then $\mathcal{R}$ is  the \dfn{reduced graph} of $G$ corresponding to the given Gallai-partition under $c$. Clearly,  $\mathcal{R}$ is isomorphic to $K_p$.  
By Theorem~\ref{Gallai},  all edges in $\mathcal{R}$ are colored by at most two colors under $c$.  One can see that any monochromatic $H$ in $\mathcal{R}$ under $c$ will result in a monochromatic $H$ in $G$ under $c$. It is not a surprise then that  Gallai-Ramsey numbers $gr_k(K_3, H)$ are related to  the classical Ramsey numbers $r_2(H)$.  Recently,  Fox,  Grinshpun and  Pach posed the following  conjecture on $gr_k(K_3, H)$ when $H$ is a complete graph.\medskip

\begin{conj}[\cite{FGP}]\label{Fox} For integers $k\ge1$ and $t\ge3$,  \\

	$gr_k(K_3, K_t) = \begin{cases}
			(r_2(K_t)-1)^{k/2} + 1 & \text{if } k \text{ is even} \\
			(t-1)  (r_2(K_t)-1)^{(k-1)/2} + 1 & \text{if } k \text{ is odd.}
			\end{cases}$
\end{conj}
\medskip

The first case of Conjecture~\ref{Fox} was verified  in 1983 due to  Chung and Graham~\cite{chgr}.  The next open case when $t=4$ was recently settled in~\cite{K4}. A simpler proof of Theorem~\ref{C3} can be found in~\cite{exponential}.

\begin{thm}[\cite{chgr}]\label{C3} For any integer $k\ge1$, 
	$gr_k(K_3, C_3) = \begin{cases}
			5^{k/2} + 1 & \text{if } k \text{ is even} \\
			2 \cdot  5^{(k-1)/2} + 1 & \text{if } k \text{ is odd.}
			\end{cases}$
\end{thm}
\medskip

  Theorem~\ref{general} below  is a result of Fujita and  Magnant~\cite{c5c6}, which provides  a lower bound for $gr_k(K_3, C_{2n+1})$. 
\begin{thm}[\cite{c5c6}]\label{general}
For  integers $k \ge 2$ and $n \ge 2$, 
$
 gr_k(K_3, C_{2n+1}) \ge n \cdot 2^k + 1.
$
\end{thm}

The exact values of $gr_k(K_3, C_5)$ for any integer $k\ge1$ were  also determined in~\cite{c5c6}.

\begin{thm}[\cite{c5c6}]\label{C5}
 $gr_k(K_3, C_{5}) = 2 \cdot 2^k + 1$ for all $k \ge 1$.
\end{thm}

Recently, Bruce and Song~\cite{DylanSong} considered the next step and determined  the exact values of $gr_k(K_3, C_7)$ for any integer $k\ge1$. 
\begin{thm}[\cite{DylanSong}]\label{C7}
 $gr_k(K_3, C_{7}) = 3 \cdot 2^k + 1$ for all $k \ge 1$.
\end{thm}

In this paper, we continue to study the Gallai-Ramsey numbers of odd cycles. 
We determine the exact values of $gr_k(K_3, C_9)$ for all  $k\ge1$. We believe the method we developed will be helpful in determining the exact values of $gr_k(K_3, C_n)$ for odd integers  $n\ge11$. We prove the following main result.

\begin{thm}\label{C9}
$gr_k(K_3, C_9) = 4 \cdot 2^k + 1$ for all $k \ge 1$.
\end{thm}

It is worth mentioning that Theorem~\ref{C9}  also provides partial evidence for the first open case of the Triple Odd Cycle Conjecture due to Bondy and Erd\H{o}s~\cite{BE}, which states that $r_3(C_n)=4n-3$ for any odd integer $n>3$. \L uczak~\cite{Luczak} showed that if $n$ is odd, then $r_3(C_n) =4n + o(n)$, as $n\rightarrow \infty$,  and Kohayakawa,
Simonovits and Skokan~\cite{TOCC} proved that the Triple Odd Cycle Conjecture holds when $n$ is sufficiently large.   We will make use of the following result of Bondy and Erd\H{o}s~\cite{BE}.

\begin{thm}[\cite{BE}]\label{2n+1}
 $r_2(C_{2n+1}) =4n+1$ for all $n\ge2$.
\end{thm}

Finally, we need to introduce more notation. For positive  integers $n, k$ and $G=K_n$, let $c$ be any $k$-edge-coloring of $G$ with color classes $E_1, \dots, E_k$. Then $c$ is  \dfn{bad}  if $G$ contains  neither a rainbow $K_3$ nor a monochromatic $C_9$ under  $c$. 
For any $E\subset E(G)$, let $G[E]$ denote the subgraph of $G$ with vertex set $V(E)$ and edge set $E$.  Let $H$ be an induced subgraph of $G$ and let  $E=E_i\cap E(H)$ for some $i\in [k]$.  Then  
$G[E]$ is an \dfn{induced matching} in $H$ if  $E$ is a matching in $H$.
For  two disjoint sets $A, B\subseteq V(G)$,  if all the edges between $A$ and $B$  in $G$ are colored the same color under $c$, say, blue,  we say that $A$ is \dfn{blue-complete} to $B$. 

\section{Proof of Theorem \ref{C9}}
By  Theorem \ref{general}, $gr_k(K_3, C_9) \ge 4 \cdot 2^k + 1$ for all $k \ge 1$. We next show that $gr_k(K_3, C_9) \le 4 \cdot 2^k + 1$ for all $k \ge 1$.  This is trivially true for $k=1$. By Theorem \ref{2n+1}, we may assume that $k\ge3$. Let $G=K_{4\cdot 2^k + 1}$ and let $c$ be any $k$-edge-coloring of $G$ such that $G$ admits no rainbow triangle.  
We next show that  $G$ contains a  monochromatic $C_9$  under the coloring $c$.\medskip
	
 Suppose that  $G$ does not  contain a  monochromatic $C_9$ under $c$.  Then $c$ is  bad. Among all complete graphs on $4\cdot 2^k+1$ vertices with a bad $k$-edge-coloring,  we choose $G$ with $k$ minimum.  We next  prove several claims.

\begin{claim}\label{3-vertex}
Let $H$ be an induced subgraph of $G$.  If there exist three distinct vertices $u, v, w \in V(G \setminus H)$ such that all edges between $\{u, v, w\}$ and $V(H)$ are colored, say blue, under $c$,  then 
\begin{enumerate}[(i)]
\item there exists $V_H \subseteq V(H)$ with $|V_H| \le 4$ such that $H \setminus V_H$ has no blue edges, and 
\item $|H| \le 4 \cdot 2^{k-1-q}+4$,  where $q\in\{0, 1, \ldots, k-1\} $ is the number of colors missing on $E(H)$ under $c$, other than blue.
\end{enumerate}
\end{claim}

\pf   To prove (i), suppose that  for any  $V_H\subseteq V(H)$ satisfying that  $H\setminus V_H$ has  no blue edges,  $|V_H| \ge5$. Then $H$ must contain three blue edges $u_1v_1, u_2v_2, u_3v_3$ such that $u_1, u_2, u_3, v_1, v_2, v_3$ are all distinct. Thus we obtain a blue $C_9$ with vertices $u, u_1, v_1, v, u_2, v_2, w, u_3, v_3$ in order, a contradiction.\medskip

By (i),  $H \setminus V_H$ has no blue edges.   By minimality of $k$, $H \setminus V_H \le 4 \cdot 2^{k-1-q}$.
Then
$|H| = |H \setminus V_H | + |V_H| \le 4 \cdot 2^{k-1-q}  + 4$.
This proves (ii).  \qed\\

Let $X_1, X_2,  \ldots, X_m$ be a maximum sequence of disjoint subsets of $V(G)$ such that for all $j \in [m]$, $1 \le |X_j| \le 2$, and all edges between $X_j$ and $V(G) \setminus \bigcup_{i \in [j]} X_i$ are colored the same color.  Let $X:= \bigcup_{j \in [m]} X_j$.   For each $x\in X$,  let $c(x)$ be the unique color on the edges between $x$ and $V(G) \less X$.   
For all $i\in[k]$, let $X^*_i :=\{x\in X: \, c(x)=\text{ color } i\}$.  Notice that $X= \bigcup_{i \in [k]} X^*_i$, and  for any $i\in[k]$, $X^*_i$ is possibly empty.

\begin{claim}\label{X*}
 For all $i\in[k]$, $|X^*_i|\le2$.
\end{claim}

\pf  Suppose not. Then $m\ge2$. When  choosing $X_1, X_2,  \ldots, X_m$, let $\ell\in[m-1]$ be the largest  index such  that   $|X_p^*\cap (X_1\cup\dots \cup X_{\ell})|\le2$ for all $p\in[k]$; and let $ j\in\{\ell+1, \ldots, m\}$ be the smallest   index such that   $3\le |X_i^*\cap (X_1\cup\dots \cup X_j)|\le4$ for some $i\in[k]$.  Such a color $i$ exists due to the assumption that the statement is not true; and such an index $ j$ exists  because $1\le |X_p|\le 2$ for all $p\in [m]$.  We may assume that color $i$ is blue.  Let  $B:=\{x\in X_1\cup\cdots\cup X_j: \, x \text{ is blue-complete to } V(G)\less X\}$, and 
  let $A :=X_1\cup \cdots\cup X_j$. By the choice of $j$,  $3\le |B|\le4$,  and  $|X_p^*\cap A|\le 2$ for any $p\in[k]\less i$. 
  Thus  $|A\less B|\le 2(k-1)$.  By Claim \ref{3-vertex} applied to any three vertices in $B$ and the induced subgraph $G\less A$,  we see that $ |G \setminus A|\le 4 \cdot 2^{k-1} +4$. Thus, 
\[
|G| = |A\less B| + |B|+ |G \setminus A| \le 2(k-1) + 4+ 4 \cdot 2^{k-1} +4 < 4 \cdot 2^k + 1
\]
for all $k \ge 3$, a contradiction.
\qed\\

By Claim \ref{X*}, $|X| \le 2k$.  Let $X'\subseteq X$  be such that for  all $i \in [k]$,  $|X'\cap X^*_i|=1$ when $X^*_i \ne \es$, and  $|X'\cap X^*_i|=0$ when $X^*_i = \es$.  Let $X'':= X\less X'$. 
Now consider a  Gallai partition $A_1, \ldots, A_p$ of $G \setminus X$ with $p\ge2$. 
We may  assume that  $1\le |A_1| \le \cdots \le  |A_s|<3 \le  |A_{s+1}| \leq \cdots \leq |A_{p}|$, where $0\le  s\le p$. Let  $\mathcal{R}$ be  the reduced graph of $G\less X$ with vertices $a_1, a_2, \dots, a_{p}$, where $a_i\in A_i$ for all $i\in[p]$. By Theorem~\ref{Gallai}, we may assume that  the edges of $\mathcal{R}$ are colored red and blue. Notice that any monochromatic $C_9$ in $\mathcal{R}$ would yield a monochromatic $C_9$ in $G$. Thus $\mathcal{R}$ has neither red nor blue $C_9$.  By Theorem~\ref{2n+1},    $p \leq 16$. Then $|A_p|\ge3$ because $|G|\ge33$.  Thus $p-s\ge1$. Let 
\[
B := \{a_i \in V(\mc{R}) \mid a_ia_p \text{ is colored blue } \} \text{ and }  R:= \{a_i \in V(\mc{R}) \mid a_ia_p \text{ is colored red } \}.
\]
Then $|B|+|R|=p-1$. 
Let $B_G:= \bigcup_{a_i \in B} A_i$, $R_G:=\bigcup_{a_j\in R} A_j$, $B_G^*:=B_G\cup\{x\in X: x \text{ is blue-complete to } V(G)\less X\}$,  and   $R_G^*:=R_G\cup\{x\in X: x \text{ is red-complete to } V(G)\less X\}$.

\begin{claim}\label{4-5 Lemma}
For any two disjoint sets $Y, Z\subseteq V(G)$ with $|Y|,|Z| \ge 4$, if  all edges between $Y$ and $Z$  are colored by the same color, say blue,  then no vertex  in $ V(G) \setminus (Y \cup Z)$ can be blue-complete to $Y \cup Z$ in $G$.  Moreover, if $|Z| \ge 5$, then $G[Z]$ has no blue edges.
\end{claim}

\pf  Suppose that there exists a vertex $x\in  V(G) \setminus (Y \cup Z)$ such that $x$ is  blue-complete to $Y \cup Z$ in $G$.
Let $Y=\{y_1, \ldots, y_{|Y|}\}$ and $Z=\{z_1, \ldots, z_{|Z|}\}$. We may further assume that $z_1z_2$  is colored blue under $c$  if  $G[Z]$ has a blue edge.  We then obtain a blue $C_9$ with  vertices $y_1, x, z_1, y_2, z_2, y_3, z_3, y_4, z_4$ in order when $|Y|,|Z| \ge 4$ or vertices  $y_1, z_1, z_2, y_2, z_3, y_3, z_4, y_4, z_5$ in order when $|Z| \ge 5$ and $G[Z]$ has a  blue edge, a contradiction.  \qed

\begin{claim}\label{B}
If $|A_p|\ge4$ and $|B|\ge5$ (resp. $|R|\ge5$), then $|B_G|\le8$ (resp. $|R_G|\le8$). 
\end{claim}
\begin{pf} Suppose $|A_p|\ge4$ and $|B|\ge5$  but $|B_G|\ge9$. By Claim \ref{4-5 Lemma},  $G[B_G]$ has no blue edges and no vertex  in $ X$ is blue-complete to $V(G)\less X$.  Thus all edges of $\mathcal{R}[B]$ are colored red in $\mathcal {R}$. Since $|B|\ge5$, we can partition $B_G$ into two subsets $B_1$ and $B_2$ such that $|B_1|\ge5$ and $|B_2|\ge4$ and $B_1$ is red-complete to $B_2$ in $G$. Then   $G[B_1]$ has no red edges and no vertex in  $ X\cup R_G$ is red-complete to $B_1$ in $G$, else we obtain a red $C_9$, a contradiction. Thus $|B|=5$, $|B_2|=4$,   $B_1=A_i$ for some $i\in\{s+1, \ldots, p-1\}$, and $R_G$ is blue-complete to $B_1$ in $G$. Then $|A_p|\ge|B_1|\ge5$ and $A_p\cup R_G$ is blue-complete to $B_1$ in $G$. By Claim \ref{4-5 Lemma}, $G[A_p\cup R_G]$ has no  blue edges. Then $G[B_1\cup X']$ has neither blue nor red edges, and $G[A_p\cup R_G\cup X'']$ has no blue edges. By minimality  of $k$, $|B_1\cup X'|\le 4\cdot 2^{k-2}$ and $|A_p\cup R_G\cup X''|\le 4\cdot 2^{k-1}$. But then, 
\[ |G|=|B_2|+|B_1\cup X'|+|A_p\cup R_G\cup X''|\le 4+4\cdot 2^{k-2}+4\cdot 2^{k-1} <4\cdot 2^{k}+1
\]
for all $k\ge3$,  a contradiction. Hence, $|B_G|\le8$. Similarly, if $|A_p|\ge4$ and  $|R|\ge5$, then $|R_G|\le8$. \qed
\end{pf}

\begin{claim}\label{p}
$p\le9$.
\end{claim}

\begin{pf} Suppose  $p\ge10$. We may assume that $|B|\ge |R|$. Then  $|B|\ge5$ because $|B|+|R|=p-1$. Thus $|B_G|\ge |B|\ge5$. We claim that $|A_p|\ge4$.  Suppose  that   $|A_p|=3$. Then $k=3$ and so $|G|=33$.
If $|A_{p-4}|=3$ or $|A_{p-8}|\ge2$, then  either $\mc{R}[\{a_{p-4}, a_{p-3}, a_{p-2}, a_{p-1}, a_{p}\}]$ has a monochromatic triangle or $C_5$,    or  $\mc{R}[\{a_{p-8}, a_{p-7}, \dots,  a_{p-1}, a_{p}\}]$ has a monochromatic $C_5$ because $r_2(C_5)=9$. In either case, we see that  $G$ has  a monochromatic $C_9$, a contradiction.  Thus  $|A_{p-4}|\le2$ and $|A_{p-8}|=1$.   Then $|A_{p-5}|= 2$, else  $|G|\le 14 +11+6< 33$.  Since $r_2(C_4)=6$, we see that $\mc{R}[\{a_{p-5}, a_{p-4}, a_{p-3}, a_{p-2}, a_{p-1}, a_{p}\}]$ has a monochromatic, say blue, $C_4$, and so $G\less X$ has a blue $C_8$ because $|A_{p-5}|=2$. Thus no vertex in $X$ is blue-complete to $G\less X$ and so $|X|\le 2(k-1)=4$. 
But then $|G|\le 12+ 8 + 8 + 4 < 33$, a contradiction.  Hence $|A_p|\ge4$, as claimed. \medskip

By Claim~\ref{4-5 Lemma},   $G[B_G]$  has  no blue edges and no vertex  in $ X$ is blue-complete to $V(G)\less X$. Thus $|X|\le 2(k-1)$ and all edges of $\mc{R}[B]$ are colored red in   $\mc{R}$. Since $|A_p|\ge 4$ and $|B|\ge5$, by Claim~\ref{B},  $|B_G|\le 8$. If $|A_p|=4$, then $|R_G|=|G|-|B_G|-|A_p|-|X|\ge 4\cdot 2^k+1-8-4-2(k-1)>5$ and 
 $|R_G\cup X'|=|G|-|B_G|-|A_p|-|X''|\ge 4\cdot 2^k+1-8-4-(k-1)> 4\cdot 2^{k-1}+1$.   By Claim~\ref{4-5 Lemma},  $G[R_G]$ has no red edges and no vertex  in $ X$ is red-complete to $V(G)\less X$.  Thus $G[R_G\cup X']$ has no red edges and so $G[R_G\cup X']$ has a monochromatic $C_9$ by the choice of $k$, a contradiction. This proves that $|A_p|\ge5$. By Claim~\ref{4-5 Lemma},  $G[A_p]$ has no blue edges.  We next claim that $|R_G|\le 4$. Suppose  $|R_G|\ge 5$. By Claim~\ref{4-5 Lemma}, neither  $G[R_G]$ nor  $G[A_p]$ has red edges and no vertex  in $ X$ is red-complete to $V(G)\less X$.  Since $G[A_p\cup X']$ has neither red nor blue edges, we see that  $|A_p\cup X'|\le 4\cdot 2^{k-2}$ by the choice of $k$.  Then  $|R_G\cup X''|= 4\cdot 2^k+1-|B_G|-|A_p\cup X'|\ge4\cdot 2^k+1-8-4\cdot 2^{k-2}\ge 4\cdot 2^{k-1}+1$. Since  $G[R_G\cup X'']$ has no red edges,  we see that $G[R_G\cup X'']$ has a monochromatic $C_9$ by minimality  of $k$, a contradiction. This proves that 
$|R_G|\le 4$. Then  $|A_p\cup X'|=|G|-|B_G|-|R_G|-|X''|\ge (4\cdot 2^k+1)-8 -4-(k-1)>4\cdot 2^{k-1}+1$.  Since $G[A_p\cup X']$ has no blue edges, by the choice of $k$, 
$G[A_p\cup X']$ has a monochromatic $C_9$, a contradiction. 
  Therefore, $p \le 9$. \qed\medskip
  \end{pf}

\begin{claim}\label{Ap}
$|A_p|\ge5$.
\end{claim}

\pf By Claim~\ref{p}, $p\le9$ and so $|A_p|\ge 4$. If $|A_p|=4$, then $|A_{p-4}|\ge3$, else  $|G|\le 16+10+2k< 4\cdot 2^k+1$ for all $k\ge3$. Then  $\mc{R}[\{a_{p-4}, a_{p-3}, a_{p-2}, a_{p-1}, a_{p}\}]$ has a blue triangle or a red  $C_5$. Then  $G$ contains   a blue or red  $C_9$, a contradiction.   \qed

\begin{claim}\label{one}
$2\le p-s\le 4$.
\end{claim}

\pf  Clearly,  $p-s\le4$, else $\mc{R}[\{a_p, a_{p-1}, a_{p-2}, a_{p-3}, a_{p-4}\}]$ has a  monochromatic   $K_3$ or $C_5$, which would  yield a blue or red  $C_9$ in $G$. 
Next suppose $p-s\le1$. Then $p-s=1$ because $p-s\ge1$.    Thus $|A_i|\le2$ for all $i\in[p-1]$.  We may assume that $|R|\le |B|$. 
We claim that  $|B|\le3$. Suppose that $|B|\ge4$.   Then $|R|\le4$ because $|B|+|R|=p-1\le8$. Thus $|R_G|\le2|R|\le8$. Then  $|B_G|\le2|B|\le 8$ when $|B|=4$. If $|B|\ge5$, then  $ |B_G|\le8$  by Claim~\ref{B}. Thus $4\le |B|\le |B_G|\le8$. By Claim \ref{4-5 Lemma}, $G[A_p]$ has no  blue edges and no vertex  in $ X$ is blue-complete to $V(G)\less X$. Thus $|X| \le 2(k-1)$ and  $G[A_p \cup X']$ has no blue edges.  By minimality of $k$,   $|A_p \cup X'| \le 4 \cdot 2^{k-1}$.  Then
\[
|R_G| =|G|-|B_G|- |A_p \cup X'|-|X''| \ge 4\cdot 2^k+1 -8-4 \cdot 2^{k-1} -(k-1)>5, 
\]
since $k \ge 3$. By Claim \ref{4-5 Lemma}, $G[A_p]$ has no  red edges and no vertex  in $ X$ is red-complete to $V(G)\less X$.  Thus $|X''|\le k-2$ and by minimality of $k$,   $|A_p \cup X'| \le 4 \cdot 2^{k-2}$. But then 
\[
|G|=|R_G| +|B_G|+|A_p \cup X'|+|X''| \le 8+8+ 4 \cdot 2^{k-2} +(k-2)< 4\cdot 2^k+1
\]
for all $k\ge3$,  a contradiction.  Thus  $|B| \le 3$, as claimed. Then $|R|\le |B|\le3$.  Thus $|B_G|, |R_G|\le2|B|\le6$.  If $|B_G|\ge4$, then by Claim \ref{4-5 Lemma}, $G[A_p]$ has no  blue edges and no vertex  in $ X$ is blue-complete to $V(G)\less X$. Thus $|X''| \le k-1$ and  $G[A_p \cup X']$ has no blue edges.  By minimality of $k$,   $|A_p \cup X'| \le 4 \cdot 2^{k-1}$.  But then 
\[
|R_G| =|G|-|B_G|- |A_p \cup X'|-|X''| \ge 4\cdot 2^k+1 -6-4 \cdot 2^{k-1} -(k-1)>6, 
\]
a contradiction. Thus $|B_G|\le3$. Similarly, $|R_G|\le3$.   
Thus  $|B^*_G|\le5$ and $|X\less B^*_G|\le 2(k-1)$. If  $|B_G^*| \ge 3$,  by Claim \ref{3-vertex} applied to any three vertices in $B_G^*$ and the induced subgraph $G[A_p]$, $|A_p| \le 4 \cdot 2^{k-1} + 4$.  But then
\[
|G| = |B^*_G| + |R_G| + |A_p| + |X\less B^*_G| \le 5 + 3 + (4 \cdot 2^{k-1} + 4) +  2(k-1) < 4 \cdot 2^k + 1
\]
for all $k \ge 3$, a contradiction.  Thus   $|B_G^*| \le 2$. Similarly, $|R_G^*|\le 2$.  
Since $p\ge2$, we see that $|B|\ge1$. Then $1\le|B|\le2$. 
By maximality of $m$, $R\ne\emptyset$,  $|B|=2$,  and $B$ is neither blue- nor red-complete to $R$ in $\mc{R}$.  Thus  $|B|=|B_G^*|=2$ (and so no vertex in $X$ is blue-complete to $V(G)\less X$), $1\le |R|\le |R_G^*| \le 2$,  $|X|\le 2(k-2)+1$, $|X'|\le k-1$ and $|X''|\le k-2$.  \medskip 

Let $Z$ be a minimal set of vertices in $A_p$ such that $G[A_p\less Z]$ has no blue edges. Then $G[(A_p\less Z)\cup X']$ has no blue edges. By minimality of $k$, $|(A_p\less Z)\cup X'|\le 4\cdot 2^{k-1}$ and so 
\[
|Z| = |G|-|(A_p\less Z)\cup X'| -|X''\cup B_G\cup R_G|\ge (4\cdot 2^k+1)-4\cdot 2^{k-1}-(k-2+4)
\ge (2\cdot 2^{k-1}+4)\ge 12,
\]  
since $k\ge3$. Let $P$ be a longest blue path in $G[A_p]$ with vertices $v_1, \dots, v_{|P|}$ in order and let $B_G=\{b_1, b_2\}$. Then $|P|\le5$ else we obtain a blue $C_9$ with vertices $b_1, v_1, \ldots, v_6, b_2, u$ in order, where $u\in A_p\less \{v_1, \ldots, v_6\}$, a contradiction. Suppose $|P|\ge4$. Then $G[A_p\less V(P)]$ has no blue path, say $P^*$,  on $7-|P|$ vertices, else we obtain a blue $C_9$ via the vertices in $B_G$, $V(P)$ and $V(P^*)$, a contradiction.  Since $|Z|\ge12$, we see that $|P|=4$ and all blue edges in $G[A_p\less V(P)]$ induce a blue matching. Notice that $Z\less V(P)$ contains exactly one vertex from each blue edge in $G[A_p\less V(P)]$ by the choice of  $Z$. Thus $G[Z\less V(P-v_1)]$ has no blue edges and  $|Z\less V(P-v_1)|\ge (2\cdot 2^{k-1}+4)-(|P|-1)\ge 2\cdot 2^{k-1}+1$. By Theorem~\ref{C5}, $G[Z\less V(P-v_1)]$ has a monochromatic, say green (possibly red), $C_5$. Let    $u_1, u_2, u_3,  u_4, u_5$ be the vertices of the $C_5$ in order.  Let $u_1u_1', \ldots, u_5u_5'$ be a blue matching in $G[A_p]$. This is possible because all blue edges in $G[A_p\less V(P)]$ induce a blue matching. Since $G$ has no rainbow triangles under the coloring $c$, we see that for any $i\in\{1,3\}$, $\{u_i, u_i'\}$ is green-complete to $\{u_{i+1}, u_{i+1}'\}$. Thus we obtain a green $C_9$ with vertices $u_1, u_2', u_1', u_2, u_3, u_4', u_3', u_4, u_5$ in order, a contradiction. This proves that $|P|\le3$. \medskip

\begin{figure}[htbp]
\begin{center}
\includegraphics[scale=0.8]{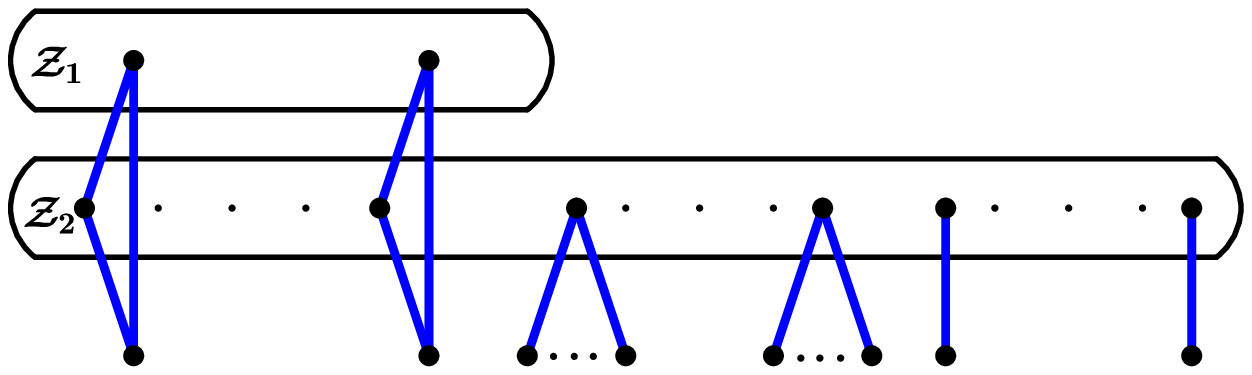}
\end{center}
\caption{Partition of $V(H)$}\label{Z1Z2}
\end{figure}

Let $H$ be the subgraph of $G[A_p]$ induced by all blue edges in $G[A_p]$. Since $|P|\le3$, we see that each  component of $H$ is isomorphic to a triangle,  a star, or a $K_2$. Let $\mc{Z}_1$ denote the set consisting of one vertex  from each $K_3$ in $H$, and   let $\mc{Z}_2$ be the set constructed by selecting: one vertex from each  $K_3\less \mc{Z}_1$   in $H$; the center vertex in each   star  in $H$; and one vertex in each $K_2$ component  in $H$, as shown in Figure~\ref{Z1Z2}.   
Finally, let $\mc{Z}_3:=A_p\less (\mc{Z}_1\cup \mc{Z}_2)$.  
Since no vertex in $X$ is blue-complete to $V(G)\less X$, neither $G[\mc{Z}_2\cup X']$ nor $G[\mc{Z}_3\cup X'']$  has blue edges.  Suppose $|\mc{Z}_2\cup X'| \ge 2 \cdot 2^{k-1} + 1$. By Theorem~\ref{C5},  $G[\mc{Z}_2\cup X']$ has  a monochromatic, say green (possibly red), $C_5$ with vertices, say $b_1, b_2, b_3, b_4, b_5$,  in order.  By the choice of $X'$, $|X'\cap V(C_5)|\le1$. We may assume that $b_1, b_2, b_3, b_4\notin X'$. Let $b_1', b_2', b_3', b_4'\in V(H)$ be such that $b_1b_1', b_2b_2', b_3b_3', b_4b_4'\in E(H)$. By the choice of $\mc{Z}_2$, $b_1', b_2', b_3', b_4'$ are all distinct.  Since $G$ has no rainbow triangle under $c$,  $\{b_i, b_i'\}$ is green-complete to $\{b_{i+1}, b_{i+1}'\}$ in $G$ for all $i\in\{1,3\}$.  We then obtain a green $C_9$ in $G$ with vertices  $b_1, b_2', b_1', b_2, b_3, b_4', b_3', b_4, b_5$ in order, a contradiction. Thus 
 $|\mc{Z}_2\cup X'| \le 2 \cdot 2^{k-1}$. Since $G[\mc{Z}_3\cup X'']$ has no blue edges, by minimality of $k$, $|\mc{Z}_3\cup X''|\le 4\cdot 2^{k-1}$.  Hence, 
\begin{equation}\label{Z2}
\begin{split}
|\mc{Z}_1| &= |G| - |\mc{Z}_2 \cup X'| - |\mc{Z}_3 \cup X''| - |B_G\cup R_G|\\
	&\ge (4 \cdot 2^k + 1) - (2 \cdot 2^{k-1}) - (4 \cdot 2^{k-1}) -  |B_G\cup R_G|\\
	&= 2^k + 1 - |B_G \cup R_G|.
\end{split}
\end{equation}
Let $3 \mc{Z}_1$ denote the  set of vertices  of all $K_3$'s in $H$.  By \eqref{Z2}, 
\begin{equation}\label{3Z2}
\begin{split}
|3\mc{Z}_1| &\ge 3\cdot 2^k + 3 - 3|B_G\cup R_G| = (2 \cdot 2^k + 1) + (2^k + 2 - 3|B_G \cup R_G|).
\end{split}
\end{equation}
If $k \ge 4$ or $|B_G \cup R_G|\le 3$, then by \eqref{3Z2}, $ |3\mc{Z}_1| \ge 2 \cdot 2^k + 1$.  By Theorem~\ref{C5},  $G[3\mc{Z}_1]$ has  a monochromatic, say green (possibly red), $C_5$ with vertices, say $b_1, b_2, b_3, b_4, b_5$,  in order.  Since $H[3\mc{Z}_1]$ consists of disjoint copies of $K_3$'s, we may assume that $N_H(b_3)\cap V(C_5)=\emptyset$.  Let $b_3, b_3', b_3''\in V(H)$ be the vertices of the $K_3$ in $H$ containing $b_3$.  Let $b_1',  b_2', b_4', b_5'\in V(H\less \{b_3, b_3', b_3''\})$  be  such that  $b_1b_1', b_2b_2', b_4b_4',b_5b_5'\in E(H)$. Note that $b_1',  b_2', b_4', b_5'$ are not necessarily distinct. \medskip

\begin{figure}[htbp]
\centering
\subfigure[][\label{b2b4}]{
\hfill\includegraphics[scale=1.2]{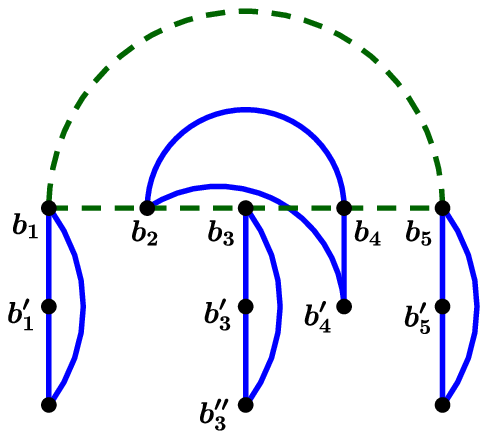}
}
\hskip 1cm
\subfigure[][\label{b2nb4}]{
\hfill\includegraphics[scale=1.2]{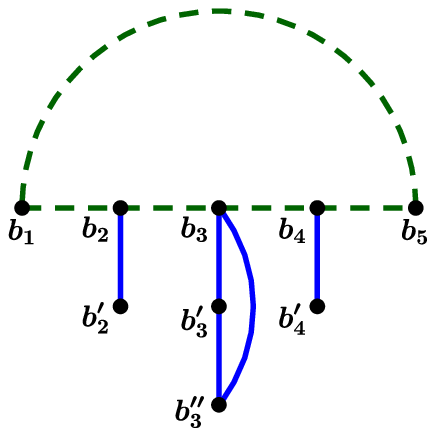}
}
\caption{When $G[3\mc{Z}_1]$ has a green $C_5$.}

\end{figure}

Suppose $b_2'= b_4'$. Then $b_1,b_3, b_4, b_5$ each belong to different $K_3$'s in $H$, as shown in Figure~\ref{b2b4}.  Since $G$ has no rainbow triangle under $c$,  $\{b_1, b_1'\}$  is green-complete to $\{b_{5}, b_{5}'\}$, and $\{b_3, b_3'\}$ is green-complete to  $\{b_4, b_4'\}$ in $G$.  We then obtain a green $C_9$ in $G$ with vertices  $b_1, b_2,  b_3, b_4', b_3', b_4, b_5', b_1', b_5$ in order, a contradiction. Thus $b_2'\ne b_4'$. Then $b_2', b_4', b_3', b_3''$ are all distinct and $b_2,b_3, b_4$ each belong to  different $K_3$'s in $H$, as shown in Figure~\ref{b2nb4}.  It  is possible that $H[\{b_1, b_4, b_4'\}]=K_3$  or $H[\{b_2, b_2', b_5\}]=K_3$.    Since $G$ has no rainbow triangle under $c$,  $\{b_3, b_3', b_3''\}$ is green-complete to $\{b_2, b_2', b_4, b_4'\}$ in $G$.  We then obtain a green $C_9$ in $G$ with vertices  $b_1, b_2,  b_3', b_2', b_3, b_4', b_3{''}, b_4, b_5$ in order, a contradiction. Thus, $k = 3$ and $|B_G \cup R_G| = 4$. 
Let the third color of $c$ be green. Thus all edges of $G[\mc{Z}_1] $ are colored red or green.  By \eqref{Z2}, $|\mc{Z}_1| \ge 8 + 1 - 4 = 5$.  Let $b_1, b_2, b_3\in \mc{Z}_1$ be distinct and let $R_G=\{u_1, u_2\}$. Let $b_i, b_i', b_i''$ be the vertices of the $K_3$ in $H$ containing $b_i$ for each $i\in[3]$. If $b_1b_2$ is colored red in $G$, then $\{b_1, b_1', b_1''\}$ must be red-complete to $\{b_2, b_2', b_2''\}$ in $G$ because $G$ has no rainbow triangle. Then we obtain a red $C_9$ in $G$ with vertices  $u_1, b_1, b_2',  b_1', b_2'', b_1'', b_2, u_2, b_3$ in order, a contradiction. Thus all edges of  $G[\mc{Z}_1] $ are colored green. Since $G$ has no rainbow triangle, $\{b_1, b_1', b_1''\}$ must be green-complete to $\{b_2, b_2', b_2'', b_3, b_3', b_3''\}$ and $\{b_2, b_2', b_2''\}$ must be green-complete to $\{ b_3, b_3', b_3''\}$.    We   obtain a green  $C_9 $ with vertices $ b_1, b_2, b_3, b_1',  b_2',  b_3',  b_1'',  b_2'',  b_3''$ in order,  a contradiction.\qed

\begin{claim}\label{A_{p-2}}
$|A_{p-2}|\le3$.
\end{claim}

\begin{pf}  Suppose   $|A_{p-2}|\ge4$. Then $4 \le |A_{p-2}| \le |A_{p-1}| \le |A_p|$ and so  $\mc{R}[\{a_{p-2}, a_{p-1},a_p\}]$ is not a monochromatic triangle in $\mc{R}$.   Let $B_1$, $B_2$, $B_3$ be a permutation of $A_{p-2}$, $A_{p-1}$, $A_p$ such that $B_2$ is, say blue-complete,  to $B_1 \cup B_3$ in $G$. Then $B_1$ must be  red-complete to $B_3$ in $G$. We may assume that  $|B_1|\ge|B_3|$.  By Claim~\ref{4-5 Lemma},    no vertex in $X$ is blue- or red-complete to $V(G)\less X$.     Let $A:=V(G)\less (B_1\cup B_2\cup B_3\cup X)$. Then by Claim~\ref{4-5 Lemma},  no vertex in $A$ is red-complete to $B_1\cup B_3$ in $G$,  and no vertex in $A$ is blue-complete to $B_1\cup B_2$ or $B_2\cup B_3$ in $G$. This implies that 
$A$ must be  red-complete to $B_2$ in $G$. We next show that  $G[A]$ has no blue edges. 
Suppose that $G[A]$ has a blue edge, say, $uv$.  Let 
\[
\begin{split}
B_1^*&:=\{b\in A \mid  b \text{ is blue-complete to   }  B_1 \text{ only} \}\\
B_2^*&:=\{b\in A\mid  b \text{ is blue-complete to both   }  B_1 \text{ and }  B_3 \}\\ 
B_3^*&:=\{b\in A\mid  b \text{ is blue-complete to   }  B_3 \text{ only} \}.\\ 
\end{split}
\]
Then $A=B_1^*\cup B_2^*\cup B_3^* $.  Notice that $B_1^*,  B_2^*, B_3^* $ are possibly empty and pairwise disjoint.   
 Let $b_1, b_2, b_3\in B_1$, $ b_4,b_5, b_6\in B_2$, and $b_7\in B_3$.  If $uv$ is an edge in $G[B_1^*\cup B_2^*]$, then we obtain a blue $C_9$ with vertices $b_1, u, v, b_2, b_4, b_7, b_5, b_3, b_6$  in order, a contradiction. 
Similarly,  $uv$ is not an edge in $G[B_2^*\cup B_3^*]$.  Thus $uv$ must be  an edge in $G[B_1^*\cup B_3^*]$ with one end in $B_1^*$ and the other in $B_3^*$.  We may assume that $u\in B_1^*$ and $v\in B_3^*$. Then we obtain a blue  $C_9$ with vertices $b_1, u, v, b_7, b_4, b_2, b_5, b_3, b_6$  in order, a contradiction. This proves that $G[A]$ has no blue edges.  \medskip

We next show that $|B_2\cup A\cup X'|\le 4\cdot 2^{k-1}$. If  $|B_2|\ge5$, then by Claim~\ref{4-5 Lemma}, $G[B_2]$ has no blue edges. Since $G[A]$ has no blue edges,  $A$ is red-complete to $B_2$,  and no vertex in $X$ is blue-complete to $B_2$, we see that   $G[B_2\cup A\cup X']$  has no blue edges. By minimality  of $k$,  $|B_2\cup A\cup X'|\le 4\cdot 2^{k-1}$. So we may assume that  $|B_2|=4$. Then $|B_2\cup A\cup X'|\le 4+4+(k-1)< 4\cdot 2^{k-1}$ when $|A|\le4$. So we may assume that  $|A|\ge 5$. By Claim~\ref{4-5 Lemma}, $G[A]$ has no red edges and no vertex in  $X$ is red-complete to $V(G)\less X$. Then $G[A\cup X']$ has neither blue nor red edges. By minimality of $k$, $|A\cup X'|\le 4\cdot 2^{k-2}$ and so $|B_2\cup A\cup X'|\le 4+4\cdot 2^{k-2}<4\cdot 2^{k-1}$. This proves that  $|B_2\cup A\cup X'|\le 4\cdot 2^{k-1}$. \medskip

 Since $|B_1|\ge |B_3|$ and  $|B_1|+|B_3|=|G|-|B_2\cup A\cup X'|-|X''|\ge 4\cdot 2^{k-1}+1-(k-1)\ge15$, we see that $|B_1|>4$. Note that  $|B_2|\ge4$ and  $|B_3|\ge4$. By Claim~\ref{4-5 Lemma},    $G[B_1]$  has neither red nor blue edges. Since each  vertex in $X$ is neither red- nor blue-complete to $B_1$,  $G[B_1\cup X'']$ has neither red nor blue edges. By minimality of $k$, $|B_1\cup X''|\le 4\cdot 2^{k-2}$ and so $|B_3|\le |B_1|\le  4\cdot 2^{k-2}$. But then  
 \[
 |G|=|B_2\cup A\cup X'|+|B_1\cup X''|+|B_3|\le 4\cdot 2^{k-1}+4\cdot 2^{k-2}+4\cdot 2^{k-2}=4\cdot 2^{k},
 \] 
   a contradiction.  \qed\bigskip
   \end{pf}

  By Claim~\ref{one},   $2\le p-s\le4$ and so $|A_{p-1}|\ge3$.  We may now assume that $a_pa_{p-1}$ is colored blue in $\mc{R}$. Then   $a_{p-1}\in B$ and so $A_{p-1}\subseteq B_G$. Thus $|B_G|\ge|A_{p-1}|\ge3$.  

\begin{claim}\label{R*_G}
$|R^*_G|\le8$.
\end{claim}

\begin{pf} 
Suppose that $|R^*_G|\ge9$.  By Claim~\ref{Ap}, $|A_p|\ge5$. By Claim~\ref{4-5 Lemma},  $G[R^*_G]$ has no red edges. Thus $|R^*_G|=|R_G|\ge9$ and so no vertex in $X$ is red-complete to $V(G)\less X$. By  Claim~\ref{B}, $|R|\le4$.  By Claim~\ref{A_{p-2}},  $|A_{p-2}|\le3$. Since $A_{p-1}\cap R_G=\emptyset$ and $p-s\le4$, we see that  $|A_{p-2}|=3$ and  $|R|=4$.   Then all edges in  $\mc{R}[R]$ are colored blue because  $G[R_G]$ has no red edges.  It can be easily checked that   $G[R_G]$ has a blue $C_9$, a contradiction.\qed
\end{pf}

\begin{claim}\label{A_{p-1}}
$|A_{p-1}|\le4$.
\end{claim}
\begin{pf} 
Suppose   $|A_{p-1}|\ge5$. Then $5 \le |A_{p-1}| \le |A_p|$.   
Thus  $|B_G|\ge|A_{p-1}|\ge5$.  By Claim~\ref{4-5 Lemma},  neither $G[A_p]$ nor $G[B_G]$  has  blue edges,  and   no vertex in $X$ is blue-complete to $V(G)\less X$.  Thus $|X|\le 2(k-1)$.  By the choice of $k$, $|B_G\cup X''|\le 4\cdot 2^{k-1}$ and $|A_p\cup X'|\le 4\cdot 2^{k-1}$. We claim that $G[R_G]$ has blue edges. Suppose that $G[R_G]$ has no blue edges. Then  $G[A_p\cup R_G\cup X']$ has no blue edges. By the choice of $k$,  $|A_p\cup R_G\cup X'|\le 4\cdot 2^{k-1}$. But then $|B_G\cup X''|=|G|-|A_p\cup R_G\cup X'|\ge 4\cdot 2^{k-1}+1$, a contradiction. Thus   $G[R_G]$ has blue edges, as claimed. Then  $|R_G|\ge2$.  Suppose that  $|R^*_G|\ge4$.  By Claim~\ref{4-5 Lemma},   $G[A_p]$   has  no red edges.  By the choice of $k$, $|A_p\cup (X'\less R^*_G)|\le 4\cdot 2^{k-2}$.  By Claim~\ref{R*_G},   $|R^*_G|\le8$.    But then
 \[  |G|=|A_p\cup (X'\less R^*_G)|+|B_G\cup (X''\less R^*_G)|+|R^*_G|\le 4\cdot 2^{k-2}+ 4\cdot 2^{k-1} +8<4\cdot 2^k+1,
  \]
for all $k\ge3$, a contradiction. 
Next suppose that  $|R^*_G|=3$. By Claim~\ref{3-vertex} applied to the three vertices in $R^*_G$ and the induced subgraph $G[A_p]$,  $|A_p|\le 4\cdot 2^{k-2}+4$. But then
 \[  |G|\le |A_p|+|B_G\cup X''|+|R^*_G|+|X'\less R^*_G|\le (4\cdot 2^{k-2}+4)+ 4\cdot 2^{k-1} +3+(k-2)<4\cdot 2^k+1,
  \]
for all $k\ge3$, a contradiction. Thus  $|R^*_G|=|R_G|=2$. Then $|X''|\le k-2$.   Let $R_G=\{a, b\}$. Then $ab$ must be colored blue under $c$ because $G[R_G]$ has blue edges.  If $a$ or $b$, say $b$,   is red-complete to $B_G$ in $G$, then neither $G[A_p\cup\{a\}\cup X']$ nor $G[B_G\cup\{b\}\cup X'']$ has  blue  edges. By minimality of $k$, $|A_p\cup \{a\}\cup X'|\le 4\cdot 2^{k-1}$ and $|B_G\cup\{b\}\cup X''|\le 4\cdot 2^{k-1}$. But then  $|G|=|A_p\cup \{a\}\cup X'|+|B_G\cup\{b\}\cup X''|\le 4\cdot 2^{k-1}+ 4\cdot 2^{k-1}<4\cdot 2^k+1$ for all $k\ge3$, a contradiction. 
Let $a', b'\in B_G$ be such that $aa'$ and $bb'$ are colored blue under $c$. Then $a'=b'$, else we obtain a blue $C_9$ in $G$ with vertices $a', a, b, b', x_1, y_1, x_2, y_2, x_3$ in order, where $x_1, x_2, x_3\in A_p$ and $y_1, y_2, y_3\in B_G\less \{a', b'\}$, a contradiction.  Thus  $\{a, b\}$  is  red-complete to $B_G\less a'$ in $G$.  Then there exists $i\in[s]$ such that $A_i=\{a'\}$.   Since $G[B_G]$ has no blue edges, we see that 
 $\{a, b, a'\}$  must be  red-complete to $B_G\less a'$ in $G$.  By Claim~\ref{3-vertex} applied to the three vertices $a, b, a'$ and the induced subgraph $G[B_G\less a']$,  we see that $|B_G\less a'|\le 4\cdot 2^{k-2}+4$. But then
 \[  |G|= |A_p\cup X'|+|B_G\less a'|+|\{a, b, a'\}|+|X''|\le  4\cdot 2^{k-1} +(4\cdot 2^{k-2}+4)+3+(k-2)<4\cdot 2^k+1,
  \]
 for all $k\ge3$, a contradiction.  Hence, $|A_{p-1}|\le4$.
 \qed\bigskip
\end{pf}

 By Claim~\ref{A_{p-2}},  $|A_{p-2}|\le3$. Thus $|A_p|\ge 6$ because $|G|\ge33$ and $p\le9$.  By  Claim \ref{one},   $2\le p-s\le4$ and so $|A_{p-1}|\ge3$. By Claim~\ref{A_{p-1}},  $3\le |A_{p-1}|\le4$.   Then $|B_G\cup R_G|\le 4+6+10=20$. 
 Clearly, $|B_G|\le 12$ by Claim~\ref{B} (when $|B|\ge5$) and the fact $p-s\le4$ (when $|B|\le4$). By Claim~\ref{R*_G}, $|R_G|\le |R^*_G|\le 8$. 
 Note that  $|B_G|\ge|A_{p-1}|\ge3$.  We first consider the case when   $ |R^*_G|\ge4$.  Since $|A_p|\ge 6$, by Claim~\ref{4-5 Lemma}, $G[A_p]$ has no red  edges.  
  If  $|B_G|=3$, then by Claim~\ref{3-vertex} applied to the three vertices in $B_G$ and the induced subgraph $G[A_p]$, $|A_p|\le 4\cdot 2^{k-2}+4$. Clearly,  $|X|\le 2k$. But then 
  \[  |G|=|A_p|+|B_G|+|R_G|+|X|\le (4\cdot 2^{k-2}+4)+ 3+8+2k<4\cdot 2^k+1,
  \]
  for all  $k\ge3$,  a contradiction. Thus  $|B_G|\ge4$. By  Claim~\ref{4-5 Lemma}, $G[A_p]$ has no blue  edges and no vertex in $X$ is blue-complete to $A_p$ in $G$.  Since $G[A_p]$ has neither red nor blue  edges,  and no vertex in $X$ is  blue-complete to $A_p$ in $G$, it follows that  $|X|\le 2(k-1)$  and  $|A_p|\le 4\cdot 2^{k-2}$ by minimality of $k$. 
  But then 
  \[
  |G|=|A_p|+|X|+(|B_G|+|R_G|)\le 4\cdot 2^{k-2}+2(k-1)+20<4\cdot 2^k+1
  \]
  for all $k\ge3$, a contradiction. \medskip
  
  It remains to consider the case when $ |R^*_G|\le3$. 
  If  $|B_G|=3$, then by Claim~\ref{3-vertex} applied to  the three vertices in $B_G$ and the induced subgraph $G[A_p]$, $|A_p|\le 4\cdot 2^{k-1}+4$. But then  
  \[
  |G|=|A_p|+|B_G|+|R_G|+|X|\le (4\cdot 2^{k-1}+4)+3+3+2k <4\cdot 2^k+1, 
  \]
  for all  $k\ge3$,  a contradiction.    Thus  $|B_G|\ge4$. By  Claim~\ref{4-5 Lemma}, $G[A_p]$ has no blue  edges and no vertex in $X$ is blue-complete to $A_p$ in $G$.  Thus $|X|\le 2(k-1)$ and $|X''\less R^*_G|\le k-2$. By minimality of $k$, $|A_p\cup X'|\le 4\cdot 2^{k-1}$. 
 But then 
   \[
   |G|\le |A_p\cup X'|+|B_G|+|R^*_G|+|X''\less R^*_G|\le 4\cdot 2^{k-1}+12+3+(k-2)<4\cdot 2^k+1
   \]
for all  $k\ge3$, a contradiction. \medskip
  
  This completes the proof of Theorem~\ref{C9}.\qed


\begin{thebibliography}{99}
\bibitem{BE}
J. A. Bondy,  P. Erd\H{o}s, Ramsey numbers for cycles in graphs, J. Combin. Theory Ser. B 14 (1973) 46--54.
%
\vspace {-0.25cm}
%
\bibitem{DylanSong} D. Bruce, Z-X. Song, Gallai-Ramsey numbers of $C_7$ with multiple colors, submitted. 
%
\vspace {-0.25cm}
%
\bibitem{CEL}  K. Cameron, J. Edmonds,   L. Lov\'asz, A note on perfect graphs, Period. Math. Hungar. 17 (1986) 173--175.
%
\vspace {-0.25cm}
%
\bibitem{chgr}
F. R. K. Chung,  R. Graham, Edge-colored complete graphs with precisely colored subgraphs, Combinatorica 3 (1983) 315 - 324.
%
\vspace {-0.25cm}
%
%
%
%
%
%
%
%
%
 \bibitem{FGP}  J. Fox, A. Grinshpun,  J. Pach, The Erd\H{o}s-Hajnal conjecture for rainbow triangles,  J. Combin. Theory Ser. B 111  (2015) 75-125.
  %
\vspace {-0.25cm}
%
\bibitem{FMO}  S. Fujita, C. Magnant,  K. Ozeki, Rainbow generalizations of Ramsey theory: a survey, Graphs   Combin. 26 (2010) 1--30.
 %
\vspace {-0.25cm}
%
\bibitem{c5c6}
S. Fujita,  C. Magnant, Gallai-Ramsey numbers for cycles, Discrete Math. 311 (2011) 1247 --1254.
%
\vspace {-0.25cm}
%
\bibitem{gallai}
T. Gallai, Transitiv orientierbare Graphen, Acta Math.  Acad. Sci.  Hung. 18 (1967) 25--66.
%
\vspace {-0.25cm}
%
 \bibitem{GS}  A. Gy\'{a}rf\'{a}s,  G. S\'{a}rk\"{o}zy, Gallai colorings of non-complete graphs, Discrete Math. 310 (2010)  977--980.
%
\vspace {-0.25cm}
%
\bibitem{exponential}
A. Gy\'{a}rf\'{a}s, G.  S\'{a}rk\"{o}zy, A. Seb\H{o},  S. Selkow, Ramsey-type results for Gallai colorings, J. Graph Theory 64 (2010), 233 - 243.
%
\vspace {-0.25cm}
%
\bibitem{upbounds} M. Hall, C. Magnant, K. Ozeki,  M. Tsugaki, Improved upper bounds for Gallai-Ramsey numbers of paths and cycles, J. Graph Theory 75 (2014) 59 --74.
%
\vspace {-0.25cm}
%
\bibitem{TOCC} Y. Kohayakawa, M. Siminovits,  J. Skokan, The $3$-colored Ramsey number of odd cycles,   Electron. Notes Discrete Math. 19 (2005) 397--402.
%
\vspace {-0.25cm}
%
\bibitem{KG} J. K\"orner,  G. Simonyi, Graph pairs and their entropies: modularity problems, Combinatorica 20 (2000)  227--240.
%
\vspace {-0.25cm}
%
\bibitem{K4} H. Liu, C. Magnant, A. Saito, I. Schiermeyer,  Y.  Shi,  Gallai-Ramsey number for $K_4$, manuscript.
%
\vspace {-0.25cm}
%
\bibitem{Luczak} T. \L uczak, $R(C_n, C_n, C_n) \le (4+o(1))n$, J. Combin. Theory Ser. B 75  (1999) 174--187. 
\end{thebibliography}
\end{document}